\documentstyle[12pt, amsfonts]{amsart}
\begin{document}
\def\R{{\mathbb R}}
\def\Z{{\mathbb Z}}
\def\C{{\mathbb C}}
\newcommand{\trace}{\rm trace}
\newcommand{\Ex}{{\mathbb{E}}}
\newcommand{\Prob}{{\mathbb{P}}}
\newcommand{\E}{{\cal E}}
\newcommand{\F}{{\cal F}}
\newtheorem{df}{Definition}
\newtheorem{theorem}{Theorem}
\newtheorem{lemma}{Lemma}
\newtheorem{pr}{Proposition}
\newtheorem{co}{Corollary}
\def\n{\nu}
\def\sign{\mbox{ sign }}
\def\a{\alpha}
\def\N{{\mathbb N}}
\def\A{{\cal A}}
\def\L{{\cal L}}
\def\X{{\cal X}}
\def\F{{\cal F}}
\def\c{\bar{c}}
\def\v{\nu}
\def\d{\delta}
\def\diam{\mbox{\rm dim}}
\def\vol{\mbox{\rm Vol}}
\def\b{\beta}
\def\t{\theta}
\def\l{\lambda}
\def\e{\varepsilon}
\def\colon{{:}\;}
\def\pf{\noindent {\bf Proof :  \  }}
\def\endpf{ \begin{flushright}
$ \Box $ \\
\end{flushright}}

\title[Hyperplane inequality for measures]{A hyperplane inequality for measures of convex bodies in $\R^n,\ n\le 4$}

\author{Alexander Koldobsky}

\address{Department of Mathematics\\ 
University of Missouri\\
Columbia, MO 65211}

\email{koldobskiya@@missouri.edu}

\begin{abstract}  Let $2\le n \le 4.$ We show that for an arbitrary measure $\mu$ with positive continuous density 
in $\R^n$ and any origin-symmetric convex body $K$ in $\R^n,$
$$\mu(K) \le \frac n{n-1}\max_{\xi \in S^{n-1}} \mu(K\cap \xi^\bot)\ \vol_n(K)^{1/n},$$
where $\xi^\bot$ is the central hyperplane in $\R^n$ perpendicular to $\xi.$  This generalizes
the hyperplane inequality in dimensions up to four to the setting of arbitrary measures in place 
of the volume. In order to prove this inequality,  we first establish stability in the affirmative 
case of the Busemann-Petty problem for arbitrary measures in the following sense:
if $\e>0,$ $K$ and $L$ are origin-symmetric convex bodies in $\R^n,\ n\le 4,$ and
for every $\xi\in S^{n-1}$
$$\mu(K\cap \xi^\bot) \le \mu(K\cap \xi^\bot) +\e,$$
then
$$\mu(K)\le \mu(L) + \frac {n\e}{n-1}\vol_n(K)^{1/n}.$$

\end{abstract}  
\maketitle

\section{Introduction}
The hyperplane problem of Bourgain \cite{Bo1}, \cite{Bo2} asks whether there exists 
an absolute constant $C$ so that for any origin-symmetric convex body $K$ in $\R^n$
\begin{equation} \label{hyper}
\vol_n(K)^{\frac {n-1}n} \le C \max_{\xi \in S^{n-1}} \vol_{n-1}(K\cap \xi^\bot),
\end{equation}
where  $\xi^\bot$ is the central hyperplane in $\R^n$ perpendicular to $\xi.$
The problem is still open, with the best-to-date estimate $C\sim n^{1/4}$ established
by Klartag \cite{Kl}, who slightly improved the previous estimate of Bourgain \cite{Bo3}.
We refer the reader to recent papers \cite{EK}, \cite{DP} for the history and 
current state of the hyperplane problem.

In the case where the dimension $n\le 4$, the inequality (\ref{hyper}) can be proved with 
the best possible constant (see \cite[Theorem 9.4.11]{G3}):
\begin{equation}\label{hyper-inter}
\vol_n(K)^{\frac {n-1}n} \le \frac{\left|B_2^n\right|^{\frac{n-1}n}}{\left|B_2^{n-1}\right|} 
\max_{\xi \in S^{n-1}} \vol_{n-1}(K\cap \xi^\bot),
\end{equation}
with equality when $K=B_2^n$ is the Euclidean ball.  Here $|B_2^n|= \pi^{n/2}/\Gamma(1+n/2)$
is the volume of $B_2^n.$ Note that the constant $C$ in (\ref{hyper-inter})
is less than 1; see (\ref{less1}). 

Inequality (\ref{hyper-inter}) follows from the affirmative answer 
to the Busemann-Petty problem in dimensions up to four.
The Busemann-Petty problem, posed in 1956 (see \cite{BP}),
asks the following question. Suppose that $K$ and $L$ are origin symmetric
convex bodies in $\R^n$ such that for every $\xi\in S^{n-1}$
$$\vol_{n-1}(K\cap \xi^\bot) \le \vol_{n-1}(L\cap \xi^\bot).$$
Does it follow that
$$\vol_n(K) \le \vol_n(L)?$$
The answer is affirmative if $n\le 4$ and negative if $n\ge 5.$
The solution was completed in the end of the 90's as the result of
a sequence of papers \cite{LR}, \cite{Ba}, \cite{Gi}, \cite{Bo4}, 
\cite{L}, \cite{Pa}, \cite{G1}, \cite{G2}, \cite{Z1}, \cite{Z2}, \cite{K1}, \cite{K2}, \cite{Z3},
\cite{GKS} ; see \cite[p. 3]{K3} or \cite[p. 343]{G3} for the history of the solution.
Applying the affirmative part of the solution to the case where $L=B_2^n,$ one
immediately gets (\ref{hyper-inter}). 

In this article we prove that inequality (\ref{hyper}) holds in dimensions up to four
with arbitrary measure in place of the volume. Let $f$ be an even continuous positive 
function on $\R^n,$ and denote by $\mu$ the measure on $\R^n$ with density $f$. 
For every closed bounded set $B\subset \R^n$ or $B\subset \xi^\bot$ define
$$\mu(B)=\int\limits_B f(x)\ dx.$$  Our extension of (\ref{hyper-inter}) is as follows:

\begin{theorem}\label{main} If $2\le n \le 4$ and $K$  
is an origin-symmetric convex body in $\R^n,$ then 
\begin{equation} \label{arbmeas}
\mu(K) \le \frac n{n-1}\max_{\xi \in S^{n-1}} \mu(K\cap \xi^\bot)\ \vol_n(K)^{1/n}.
\end{equation}
\end{theorem}

Zvavitch \cite{Zv} found a remarkable generalization of the Busemann-Petty problem
to arbitrary measures, namely, one can replace the volume by any measure with positive 
continuous density in $\R^n.$ In particular, if $n\le 4,$ then for any 
convex origin-symmetric bodies $K$ and $L$ in $\R^n$ the inequalities
$$\mu(K\cap \xi^\bot) \le \mu(L\cap \xi^\bot), \qquad \forall \xi\in S^{n-1}$$
imply
$$\mu(K)\le \mu(L). $$
Zvavitch also proved that this is generally not true if $n\ge 5,$ making the answer similar to
the original Busemann-Petty problem. 

By analogy with the volume case, one would expect that this result immediately implies (\ref{arbmeas}).
The argument, however, does not work in this setting, because the measure $\mu$ of sections of the 
Euclidean ball does not have to be a constant. Instead, to prove (\ref{arbmeas}) we establish stability 
in the affirmative part of Zvavitch's result in the following sense:

\begin{theorem} \label{stab} Let  $f$ be an even positive continuous function on $\R^n,\ 2\le n \le 4,$
$\mu$ is the measure with density $f,$ $K$ and $L$ are origin-symmetric convex
bodies in $\R^n,$ and $\e>0.$ Suppose that for every $\xi\in S^{n-1},$
\begin{equation} \label{cond}
\mu(K\cap \xi^\bot) \le \mu(L\cap \xi^\bot) +\e.
\end{equation}
Then
\begin{equation}\label{concl}
\mu(K)\le \mu(L) + \frac {n\e}{n-1}\vol_n(K)^{1/n}.
\end{equation}
\end{theorem}

Interchanging $K$ and $L,$ we get
\begin{co} \label{ineq-meas}Under the conditions of Theorem \ref{stab}, we have
$$\left|\mu(K) - \mu(L)\right| $$
\begin{equation} \label{measure}
\le \frac n{n-1}\max_{\xi \in S^{n-1}}\left|\mu(K\cap \xi^\bot) - \mu(L\cap \xi^\bot)\right| 
 \max\left(\vol_n(K)^{\frac 1n}, \vol_n(L)^{\frac1n}\right).
\end{equation}
\end{co}

Now to prove Theorem \ref{main} simply put $L=\emptyset$ in Corollary \ref{ineq-meas}.
\smallbreak
It remains to prove Theorem \ref{stab}. Note that stability in the original Busemann-Petty problem
(for the volume) was established in \cite{K5}. We discuss the relation between different stability 
estimates in the end of the paper.

\section{Preliminaries}

We use the techniques of the Fourier approach
to sections of convex bodies; see \cite{K3} and \cite{KY} for details. 
As usual, we denote by ${\cal{S}}(\R^n)$
the Schwartz space of rapidly decreasing infinitely differentiable
functions (test functions) in $\R^n,$ and
${\cal{S}}^{'}(\R^n)$ is
the space of distributions over ${\cal{S}}(\R^n).$

Suppose that $f$ is a locally 
integrable complex-valued function on $\R^n$ with {\it power growth at infinity}, 
i.e. there exists a number $ \beta>0$ so that 
$$\lim_{|x|_2\to \infty} \frac{f(x)} {|x|_2^\beta}=0,$$
where $|\cdot|_2$ is the Euclidean norm in $\R^n.$
Then $f$ represents a distribution acting by integration: 
for every $\phi\in {\mathcal S},$ 
$$\langle f, \phi \rangle = \int_{\R^n} f(x) \phi(x)\ dx.$$

The Fourier transform of a
distribution $f$ is defined by $\langle\hat{f}, \phi\rangle= \langle f, \hat{\phi} \rangle$ for
every test function $\phi.$ 

A distribution $f$ is called even homogeneous of degree $p\in \R$  
$$\langle f(x), \phi(x/\alpha) \rangle = |\alpha|^{n+p} 
\langle f,\phi \rangle$$
for every test function $\phi$ and every 
$\alpha\in \R,\ \alpha\neq 0.$  The Fourier transform of an even
homogeneous distribution of degree $p$ is an even homogeneous
distribution of degree $-n-p.$ 

We say that a distribution is  {\it positive definite} if its Fourier transform is a positive distribution in
the sense that $\langle \hat{f},\phi \rangle \ge 0$ for every non-negative test function $\phi.$
L. Schwartz's generalization of Bochner's theorem (see, for example, 
\cite[p.152]{GV}) states that a distribution is positive definite
if and only if it is the Fourier transform of 
a tempered measure on $\R^n$. Recall that \index{tempered measure}
a (non-negative, not necessarily finite) measure $\mu$ is  
called tempered if 
$$\int_{\R^n} (1+|x|_2)^{-\beta}\ d\mu(x)< \infty$$
for some $\beta >0.$ 

For an origin-symmetric convex body $K$ in $\R^n$ we denote by
$$\|x\|_K = \min\{a\ge 0:\ x\in aK\}, \qquad x\in \R^n$$
the norm in $\R^n$ generated by $K.$ Our definition of a convex body
assumes that the origin is an interior point of $K.$ If  $0<p<n,$
then $\|\cdot\|_K^{-p}$  is a locally integrable function on $\R^n$ and represents an even 
homogeneous of degree $-p$ distribution. If $\|\cdot\|_K^{-p}$ represents a positive definite
distribution for some $p\in (0,n),$ then its Fourier transform is a tempered measure which 
is at the same time a homogeneous distribution of degree $-n+p.$ One can express such 
a measure in polar coordinates:

\begin{pr} \label{posdef}  (\cite[Corollary 2.26]{K3}) 
Let $K$ be an origin-symmetric convex body in $\R^n$ and
$p\in (0,n).$ The function $\|\cdot\|_K^{-p}$ represents a 
positive definite distribution on $\R^n$ if and only if there 
exists a finite Borel measure $\mu_0$ on $S^{n-1}$ so that 
for every even test function $\phi,$
$$
\int_{\R^n} \|x\|_K^{-p} \phi(x)\ dx = \int_{S^{n-1}} \left(
\int_0^\infty t^{p-1} \hat\phi(t\xi) dt \right) d\mu_0(\xi).
$$
\end{pr}

It was proved in \cite{GKS} (see \cite[Corollary 4.9]{K3}) that 
\begin{pr} If $2\le n \le 4$ and $K$ is any origin-symmetric convex body in $\R^n,$ 
then the function $\|\cdot\|_K^{-1}$ represents a positive definite distribution.
\end{pr}

For any continuous function $f$ on the sphere $S^{n-1}$ denote by $f\cdot r^{-n+1}$ the 
extension of $f$ to an even homogeneous function of degree $-n+1$ on the whole $\R^n.$
It was proved in \cite[Lemma 3.7]{K3} that the Fourier transform of $f\cdot r^{-n+1}$ 
is equal to another continuous function $g$ on $S^{n-1}$ extended to an even homogeneous 
of degree $-1$ function $g\cdot r^{-1}$ on the whole $\R^n$
(in fact, $g$ is the spherical Radon transform of $f$, up to a constant).
This is why we can remove smoothness conditions in the Parseval formula on 
the sphere \cite[Corollary 3.23]{K3} and formulate it as follows:

\begin{pr}  \label{parseval} Let $K$ be an origin-symmetric convex body in $\R^n.$ Suppose that 
$\|\cdot\|_K^{-1}$ is a positive definite distribution, and let $\mu_0$ be the finite Borel measure on $S^{n-1}$ 
that corresponds to $\|\cdot\|_K^{-1}$ by Proposition \ref{posdef}. Then for any even continuous function $f$
on $S^{n-1},$
\begin{equation} 
\int_{S^{n-1}} (f\cdot r^{-n+1})^\wedge(\theta)\ d\mu_0(\theta) =
\int_{S^{n-1}}  \|\theta\|_K^{-1} f(\theta)\ d\theta.
\end{equation}
\end{pr}
Finally, we need a formula from \cite{Zv}, expressing the measure of a section in terms 
of the Fourier transform. This formula generalizes the corresponding result for volume;
see \cite{K4}.
\begin{pr}\label{tm:f} {\bf (\cite{Zv})}
Let $K$ be an  origin-symmetric star body in $\R^n$, then, for every $\xi\in S^{n-1},$ 
$$\mu(K\cap \xi^\bot)=\frac{1}{\pi}\left(|x|_2^{-n+1}
\int\limits_{0}^{|x|_2/\|x\|_K} t^{n-2} 
f\left(\frac{tx}{|x|_2}\right)dt\right)^{\wedge}(\xi),$$
where the Fourier transform of the function of $x\in \R^n$ in the right-hand
side is a continuous homogeneous of degree $-1$ function on $\R^n\setminus \{0\}.$ 
\end{pr}

\section{Stability}
We need a simple property of the $\Gamma$-function.
\begin{lemma} \label{gammafunction} For any $n\in \N,$ 
$$ \frac{\Gamma(\frac{n-1}2)}{\left(\Gamma(\frac n2)\right)^{\frac{n-1}n}} \le \frac{n^{\frac{n-1}n}2^{\frac 1n}}{n-1}.$$
\end{lemma}
\pf By log-convexity of the $\Gamma$-function (see \cite[p.30]{K3}),
$$\frac{\log\left(\Gamma(\frac n2+1)\right)-\log\left(\Gamma(1)\right)}{n} \ge 
\frac {\log\left(\Gamma(\frac {n-1}2+1)\right)-\log\left(\Gamma(1)\right)}{n-1},$$
so
\begin{equation} \label{less1}
\Gamma\left(\frac{n+1}2\right) \le \left(\Gamma\left(\frac n2+1\right)\right)^{\frac{n-1}n}.
\end{equation}
Now apply the property $\Gamma(x+1)=x\Gamma(x)$ of the $\Gamma$-function. 
\endpf
The following elementary fact was used by Zvavitch \cite{Zv} in his generalization 
of the Busemann-Petty problem.
\begin{lemma}\label{ell} 
Let $a,b>0$ and let $\alpha$  be a non-negative 
function on $(0, \max\{a,b\}]$ 
so that the integrals below converge. Then
\begin{equation}\label{e:ell}
\int\limits_0^a t^{n-1} \alpha(t)dt - a\int\limits_{0}^a t^{n-2} 
\alpha(t)dt\le \int\limits_0^b t^{n-1} \alpha(t)dt -
a\int\limits_{0}^b t^{n-2} \alpha(t) dt.
\end{equation}
\end{lemma}
\pf The  inequality (\ref{e:ell}) is equivalent to
$$
a\int\limits_{a}^b t^{n-2} \alpha(t)dt
\le \int\limits_a^b t^{n-1} \alpha(t)dt.
$$
Note that the latter inequality also 
holds in the case $a\ge b$.
\endpf

The measure of a body can be expressed in polar coordinates as follows: 
\begin{equation} \label{polar-measure}
\mu(K) = \int_K f(u)\ du = \int\limits_{S^{n-1}}\left(\int\limits_0^{\|x\|^{-1}_K} t^{n-1} f(tx)dt\right) dx.
\end{equation}
In particular, if $f=1$ we get the polar formula for volume:
\begin{equation} \label{polar-volume}
n\vol_n(K)= \int_{S^{n-1}} \|x\|_K^{-n} dx .
\end{equation}
\smallbreak
We are ready to prove Theorem \ref{stab}.
\smallbreak
{\bf Proof of Theorem \ref{stab}.}  First, we rewrite the condition (\ref{cond}) using Proposition \ref{tm:f}:
$$
\left(|x|_2^{-n+1}\int\limits_0^{\frac{|x|_2}{\|x\|_K}} t^{n-2} 
f\left(\frac{tx}{|x|_2}\right)dt \right)^\wedge(\xi)$$$$
\le
\left(|x|_2^{-n+1}\int\limits_0^{\frac{|x|_2}{\|x\|_L}} t^{n-2} 
f\left(\frac{tx}{|x|_2}\right)dt \right)^\wedge(\xi) + \pi \e
$$
for each $\xi\in S^{n-1}.$

We integrate the latter inequality over $S^{n-1}$ with respect to the
measure $\mu_0$ corresponding to the positive definite homogeneous of degree
$-1$ distribution $\|\cdot\|_K^{-1}$ by Proposition \ref{posdef}:
$$\int\limits_{S^{n-1}} \left(|x|_2^{-n+1}\int\limits_0^{\frac{|x|_2}{\|x\|_K}} t^{n-2} 
f\left(\frac{tx}{|x|_2}\right)dt \right)^\wedge(\xi) \,d \mu_0(\xi)$$ 
$$\le\int\limits_{S^{n-1}}\left(|x|_2^{-n+1}\int\limits_0^{\frac{|x|_2}{\|x\|_L}} t^{n-2} 
f\left(\frac{tx}{|x|_2}\right)dt \right)^\wedge(\xi)\,d\mu_0(\xi) + \pi \e \int_{S^{n-1}} d\mu_0(\xi),
$$
and now apply the spherical Parseval formula, Proposition \ref{parseval}:
$$\int\limits_{S^{n-1}}\|x\|_K^{-1}\left(\int\limits_{0}^{\|x\|^{-1}_K} t^{n-2}
 f(tx)dt\right)\ dx$$
\begin{equation}\label{zvavitch1}
\le\int\limits_{S^{n-1}}\|x\|_K^{-1}\left(\int\limits_{0}^{\|x\|^{-1}_L} t^{n-2} f(tx) 
 dt\right)\ dx  + \pi \e \int_{S^{n-1}} d\mu_0(\xi).
\end{equation}

By Lemma \ref{ell} with $a= \|x\|_K^{-1}$,  $b= \|x\|_L^{-1}$ and
 $\alpha(t)=f(tx),$
\begin{align}
\int\limits_0^{\|x\|_K^{-1}}&t^{n-1} f(tx)\ dt - \|x\|_K^{-1}
 \int\limits_{0}^{\|x\|_K^{-1}} t^{n-2} f(tx)\ dt\nonumber\\
\le& \int\limits_0^{\|x\|_L^{-1}} t^{n-1} f(tx)dt -
\|x\|_K^{-1}
 \int\limits_{0}^{\|x\|_L^{-1}}t^{n-2} f(tx) dt, 
\,\,\,\, \forall x\in S^{n-1}.\nonumber
\end{align}
Integrating over $S^{n-1}$ we get
\begin{align}\label{zvavitch2}
&\int\limits_{S^{n-1}}\left(\int\limits_0^{\|x\|_K^{-1}}t^{n-1} f(tx)dt\right)\ dx  
- \int\limits_{S^{n-1}}\|x\|_K^{-1}\left( \int\limits_{0}^{\|x\|_K^{-1}}t^{n-2}
 f(tx)dt\right)\ dx\\
&\le \int\limits_{S^{n-1}}\left(\int\limits_0^{\|x\|_L^{-1}}t^{n-1} f(tx)dt\right)\ dx-
\int\limits_{S^{n-1}} \|x\|_K^{-1} \left( \int\limits_{0}^{\|x\|_L^{-1}}t^{n-2} 
f(tx)  dt\right) \ dx. \nonumber
\end{align}
Adding inequalities (\ref{zvavitch1}) and (\ref{zvavitch2}) we get
$$
\int\limits_{S^{n-1}}\left(\int\limits_0^{\|x\|^{-1}_K} t^{n-1} f(tx)dt\right) dx$$
$$
\le \int\limits_{S^{n-1}}\left(\int\limits_0^{\|x\|^{-1}_L} t^{n-1} f(tx)dt \right)dx + \pi \e \int_{S^{n-1}} d\mu_0(\xi),
$$
and, by the polar formula (\ref{polar-measure}), the latter can be written as
$$\mu(K) \le \mu(L) + \pi \e \int_{S^{n-1}} d\mu_0(\xi).$$

It remains to estimate the integral in the right-hand side of the latter inequality. For this we use the formula for the Fourier transform 
(in the sense of distributions; see \cite[p.194]{GS})
$$\left(|x|_2^{-n+1}\right)^\wedge(\xi) = \frac{2\pi^{\frac {n+1}{2}}}
{\Gamma(\frac{n-1}2)} |\xi|_2^{-1}.$$
Again using Parseval's formula, Proposition \ref{parseval},
$$\pi \e \int_{S^{n-1}} d\mu_0(\xi) = \frac{\pi \e \Gamma(\frac{n-1}2)}{2\pi^{\frac {n+1}{2}}} \int_{S^{n-1}} \left(|\cdot|_2^{-n+1}\right)^\wedge(\xi)
d\mu_0(\xi)$$
$$= \frac{\pi \e \Gamma(\frac{n-1}2)}{2\pi^{\frac {n+1}{2}}} \int_{S^{n-1}} \|x\|_K^{-1} dx
\le \frac{\pi \e \Gamma(\frac{n-1}2)}{2\pi^{\frac {n+1}{2}}} \left(\int_{S^{n-1}} \|x\|_K^{-n} dx\right)^{1/n} \left|S^{n-1}\right|^{\frac {n-1}n},$$
where we also used H\"older's inequality. Here
$$\left|S^{n-1}\right|= \frac{2\pi^{\frac{n}2}}{\Gamma(\frac{n}2)}$$
is the surface area of the unit sphere $S^{n-1}.$ Now use the polar formula for volume (\ref{polar-volume})
and apply Lemma \ref{gammafunction} to get the desired estimate.
\endpf
Theorem \ref{stab} does not hold true in dimensions greater than four, simply because the answer
to the Busemann-Petty problem in these dimensions is negative.

Stability in the original Busemann-Petty problem was studied in \cite{K5}, where it was shown that
if the dimension $n\le 4,$ then for any origin-symmetric convex bodies $K$ and $L$ in $\R^n$ and 
every $\e>0,$ the inequalities 
$$\vol_{n-1}(K\cap \xi^\bot) \le \vol_{n-1}(L\cap \xi^\bot)  + \e, \qquad  \forall \xi\in S^{n-1}$$
imply 
\begin{equation}\label{volume}
\vol_n(K)^{\frac{n-1}n} \le \vol_n(L)^{\frac{n-1}n} + \e.
\end{equation} 
This is stronger than what Theorem \ref{stab} provides in the case of the volume. In fact, if $\mu$
in Theorem \ref{stab} is the volume ($f\equiv 1$), then (\ref{concl}) reads as
$$\vol_n(K) \le \vol_n(L) + \frac {n\e}{n-1} \vol_n(K)^{1/n},$$
which follows from  (\ref{volume}) by the Mean Value Theorem applied to the function $h(t)=t^{n/(n-1)}.$

However, Theorem \ref{stab} works for arbitrary measures, while the approach of \cite{K5} does 
not allow this degree of generality.

\bigbreak
{\bf Acknowledgement.} The author wishes to thank
the US National Science Foundation for support through 
grants DMS-0652571 and DMS-1001234.

\end{document}